\documentclass[12pt]{article}
\usepackage{amsfonts,amsmath}
\def\cA{{\cal A}}          \def\cB{{\cal B}}          \def\cC{{\cal C}}
                    \def\cF{{\cal F}}
                    
          \def\cK{{\cal K}}          
                    
                    \def\cR{{\cal R}}

\newcommand{\ZZ}{{\mathbb Z}}

\newcommand{\sltwo}{{(\widehat{sl(2)}_{c})}}
\newcommand{\elpa}{{{\cal A}_{q,p}}}
\newcommand{\elpb}{{{\cal B}_{q,p,\lambda}}}
\newcommand{\un}{\mbox{1\hspace{-1mm}I}}
\newcommand{\id}{\mbox{id}}
\newcommand{\uq}[1]{{{\cal U}_{#1}}}

\newcommand{\dy}[2]{{{\cal D}Y_{#1}^{#2}}}

\def\cotan{\mathop{\rm cotg}\nolimits}
\def\Ad{\mathop{\rm Ad}\nolimits}

\def\QTHA{Quasitriangular Hopf Algebra (QTHA)\def\QTHA{QTHA}}
\def\QTQHA{Quasitriangular Quasi-Hopf Algebra (QTQHA)\def\QTQHA{QTQHA}}

\begin{document}

\pagestyle{empty}

\begin{center}

{\Large \textbf{On the quasi-Hopf structure of deformed double Yangians}}

\vspace{10mm}

{\large D. Arnaudon$^a$, J. Avan$^b$, L. Frappat$^a$, {\'E}. Ragoucy$^a$, M.
Rossi$^a$}

\vspace{10mm}

\emph{$^a$ Laboratoire d'Annecy-le-Vieux de Physique Th{\'e}orique LAPTH}

\emph{CNRS, UMR 5108, associ{\'e}e {\`a} l'Universit{\'e} de Savoie}

\emph{LAPP, BP 110, F-74941 Annecy-le-Vieux Cedex, France}

\vspace{7mm}

\emph{$^b$ LPTHE, CNRS, UMR 7589, Universit{\'e}s Paris VI/VII, France}

\end{center}

\vfill
\vfill

\begin{abstract}
  We construct universal twists connecting the centrally extended double
  Yangian $\dy{}{}(sl(2))_c$ with deformed double Yangians
  $\dy{r}{}(sl(2))_c$, thereby establishing the quasi-Hopf
  structures of the latter.   
\end{abstract}

\vfill
\vfill

\rightline{LAPTH-773/00}
\rightline{PAR-LPTHE 00-01}
\rightline{math.QA/0001034}
\rightline{January 2000}

\cleardoublepage

\baselineskip=16pt
\newpage
\pagestyle{plain}
\setcounter{page}{1}

\section{Introduction}
\setcounter{equation}{0}

Universal twists connecting (affine) quantum groups to (elliptic)
(dynamical) (affine) algebras have been constructed in
\cite{Bab,ABRR,JKOS}. They show in particular the quasi-Hopf structure
of elliptic and dynamical algebras. These twists transform the
universal $R$-matrix $\cR$ of the first object into the universal
$R$-matrix $\cR^\cF$ of the second one as follows:
\begin{equation}
  \cR^\cF_{12} = \cF_{21} \cR_{12} \cF_{12}^{-1} \;.
  \label{eq:RFRFuni}
\end{equation}

The double degeneracy limits of elliptic $R$-matrices, whether vertex-type
\cite{Konno,JKM,KLP} or face-type \cite{Cla} give rise to algebraic
structures which have been variously characterised as scaled elliptic
algebras \cite{JKM,KLP}, or double Yangian algebras
\cite{Cla,BL,KT}.  As 
pointed out earlier \cite{Konno,KLP}\footnote{We wish to thank
  S.~Pakuliak for clarifying this point to us.}, 
although represented by formally
identical Yang--Baxter relations $RLL = LLR$ \cite{FRT}, these two classes
of objects differ fundamentally in their structures (as is reflected
in the very different mode expansions of $L$ defining their individual
generators) and must be considered separately.

In our previous paper \cite{Cla} we have defined, in the quantum inverse
scattering or RLL formulation, various algebraic structures of double
Yangian type connected by twist-like operators, i.e. such that their
evaluated $R$-matrices were related as:
\begin{equation}
        R^F_{12} = F_{21} R_{12} F_{12}^{-1}
\label{eq:RFRF}
\end{equation}
for a particular matrix $F$.  It was conjectured that these twist-like
operators were indeed evaluation representations of universal twists
obeying a shifted cocycle condition thereby raising the relation
(\ref{eq:RFRF}) to the status of a genuine twist connection
(\ref{eq:RFRFuni}) between quasi-Hopf algebras. 

We shall be concerned here only with algebraic structures related to
the algebra $\widehat{sl(2)}_c$, and henceforth dispense with
indicating it 
explicitly: for instance $\dy{}{}$ is thus to be understood as 
$\dy{}{}\sltwo$. 

It is our purpose here to establish such connections, at the level of
universal $R$-matrices, between the double Yangian structures respectively
known as $\dy{}{}$, $\dy{r}{V6}$, $\dy{r}{V8}$ and $\dy{r}{F}$.  $\dy{}{}$
is the double Yangian defined in \cite{KT,Kh}.  $\dy{r}{V6}$ is
characterised by a scaled ``elliptic'' $R$-matrix defined in \cite{Konno},
$\dy{r}{V8}$ is characterised by a scaled $R$-matrix defined in
\cite{KLP,JKM}.  
In connection with our previous caveat,
note that these $R$-matrices are also used to describe
respectively the scaled elliptic algebras $\cA_{\hbar,0}$,
$\cA_{\hbar,\eta}$ 
\cite{JKM,KLP,Konno}.  $\dy{r}{F}$ is the deformed double Yangian obtained
by a particular limit of the dynamical $R$-matrix characterising elliptic
$\elpb$ algebra \cite{Cla}. 

A crucial ingredient for our procedure is
a linear (difference) equation obeyed by the twist. This type of
equation for twist operators was first written in \cite{BR}. It is
also present in \cite{ABRR,JKOS}. 
Our method consists in 
\textit{i)} finding a twist-like action in representation 
\textit{ii)} interpreting this representation as an infinite product
\textit{iii)} defining a linear equation obeyed by this infinite
product
\textit{iv)} promoting this linear equation for the representation to
the level of linear equation for universal twist. 
\textit{v)} The solution of this linear equation
is obtained as a infinite product as in \cite{ABRR}
which \textit{vi)} is then proved to obey the
shifted cocycle condition as in \cite{ABRR,JKOS}
and \textit{vii)} has an evaluation representation identical to the
twist-like action found in \textit{i)}.

This provides us with the universal $R$-matrix and quasi-Hopf structure of the
twisted algebras $\dy{r}{V6,V8,F}$,
thereby realising a fully consistent description of these algebraic
structures.

The universal $R$-matrix and Hopf algebra structure for $\dy{}{}$ were
described in \cite{KT,Kh}.  
We construct a universal twist between $\dy{}{}$ and
$\dy{r}{V6}$. We then show the existence of a universal coboundary
(trivial) twist, the evaluation of which realises the connection between the
evaluated $R$-matrices of $\dy{r}{V6}$ and $\dy{r}{V8}$, leading to
identification of these two as quasi-Hopf algebras.  Finally another
universal coboundary-like twist realises, when evaluated, the connection
between the $R$-matrices of $\dy{r}{V6}$ and $\dy{r}{F}$.

It follows that the three deformed structures are in fact one single
quasi-Hopf algebra described by three different choices of generators,
more precisely given in three different gauges. 

We shall denote throughout this paper 
${\cF}[{\cA}\,;{\cB}]$ the
universal or represented twist connecting $R$-matrices as 
${\cR}_{{\cB}} = {\cF}_{21}[{\cA}\,;{\cB}] \; {\cR}_{{\cA}} \; 
{\cF}_{12}^{-1}[{\cA}\,;{\cB}]$.

\section{Presentation of the double Yangians $\dy{}{}$ and $\dy{r}{}$} 
\setcounter{equation}{0}

\subsection{Double Yangian $\dy{}{}$}
The double Yangian $\dy{}{}$ is defined by the $R$-matrix
\begin{equation}
  R(\beta) = \rho(\beta) \left( \begin{array}{cccc} 
  1 & 0 & 0 & 0 \\
  0 & \displaystyle \frac{i\beta}{i\beta+\pi} & \displaystyle 
  \frac{\pi}{i\beta+\pi} & 0 \\[.3cm]
  0 & \displaystyle \frac{\pi}{i\beta+\pi} & \displaystyle 
  \frac{i\beta}{i\beta+\pi} & 0 \\
  0 & 0 & 0 & 1 \\
\end{array} \right) \;, 
\label{eq:dy}
\end{equation}
with the normalisation factor 
\begin{equation}
  \rho(\beta) = \frac{\Gamma_{1}(\frac{i\beta}{\pi} \;\vert\; 2) \; 
    \Gamma_{1}(2+\frac{i\beta}{\pi} \;\vert\; 2)}
  {\Gamma_{1}(1+\frac{i\beta}{\pi} \;\vert\; 2)^2} \;,
\end{equation}
together with the relations 
\begin{eqnarray}
  R_{12}(\beta_{1}-\beta_{2}) \, L^\pm_1(\beta_{1}) \,
  L^\pm_2(\beta_{2}) &=&  
  L^\pm_2(\beta_{2}) \, L^\pm_1(\beta_{1}) \, 
  R_{12}(\beta_{1}-\beta_{2}) \;. \\
  R_{12}(\beta_{1}-\beta_{2}-i\pi c) \, L^-_1(\beta_{1}) \,
  L^+_2(\beta_{2}) &=&  
  L^+_2(\beta_{2}) \, L^-_1(\beta_{1}) \, 
  R_{12}(\beta_{1}-\beta_{2}) \;. 
  \label{eq:rll_dy}
\end{eqnarray}
and the mode expansions
\begin{equation}
  L^+(\beta) = \sum_{k \ge 0} L^+_{k} \beta^{-k} \qquad \mbox{and} \qquad 
  L^-(\beta) = \sum_{k \le 0} L^-_{k} \beta^{-k} \;.
  \label{eq:mode}
\end{equation}
It is important to point out that $L^+$ and $L^-$ are independent.
There exists in this case a Gauss decomposition of the Lax matrices
allowing for an alternative Drinfeld presentation \cite{Kh}. 

\noindent
Indeed, $L^\pm$ are decomposed as 
\begin{equation}
  L^\pm(x) = 
  \left( 
    \begin{array}{cc} 1 & f^\pm(x^\mp) \\ 0 & 1 \\ 
    \end{array} 
  \right) 
  \left( 
    \begin{array}{cc} k_1^\pm(x) & 0 \\ 0 & k_2^\pm(x) \\ 
    \end{array} 
  \right) 
  \left( 
    \begin{array}{cc} 1 & 0 \\ e^\pm(x) & 1 \\ 
    \end{array} 
  \right) 
  \label{eq:Gaubeta}
\end{equation}
with $x^+\equiv x\equiv \frac{i\beta}{\pi}$ and $x^-\equiv x-c$.
Furthermore, $k_1^\pm(x) k_2^\pm(x-1) = 1$ and one defines
$h^\pm(x)\equiv k_2^\pm(x)^{-1} k_1^\pm(x)$.

\noindent 
The evaluation representation $\pi_x$ is then easily defined by
its action on a two-dimensional vector space by
\begin{equation}
  \pi_x(e_k) =x^k \sigma^+  \;,
  \qquad 
  \pi_x(f_k) =x^k \sigma^-  \;,
  \qquad
  \pi_x(h_k) =x^k \sigma^3  \;,
  \label{eq:pi}
\end{equation}
where 
\begin{equation}
  e^\pm(u) \equiv \pm\sum_{
    \begin{array}{c} 
      {\scriptstyle k\geq 0} \\[-.9ex] {\scriptstyle k<0 } \\
    \end{array}
    }
  e_k u^{-k-1}\,,
  \qquad
  f^\pm(u) \equiv \pm\sum_{
    \begin{array}{c} 
      {\scriptstyle k\geq 0} \\[-.9ex] {\scriptstyle k<0 } \\
    \end{array}
    }
  f_k u^{-k-1}\,,
  \qquad
  h^\pm(u) \equiv 1 \pm\sum_{
    \begin{array}{c} 
      {\scriptstyle k\geq 0} \\[-.9ex] {\scriptstyle k<0 } \\
    \end{array}
    }
  h_k u^{-k-1}\,.
  \label{eq:modesefh}
\end{equation}

\subsection{Deformed double Yangian $\dy{r}{V6}$}
The $R$-matrix of the deformed double Yangian $\dy{r}{V6}$ is related to
the two-body $S$ matrix of the 
sine--Gordon theory $S_{SG}(\beta,r)$ and is given by
\begin{equation}
  R_{V6}(\beta,r) = \cotan(\frac{i\beta}{2}) S_{SG}(\beta,r)
  =
  \rho_r(\beta) 
  \left( 
    \begin{array}{cccc} 
      1 & 0 & 0 & 0 \\[1mm]
      0 & \displaystyle
      \frac{\sin\frac{i\beta}{r}}{\sin\frac{\pi+i\beta}{r}} &  
      \displaystyle \frac{\sin\frac{\pi}{r}}
      {\sin\frac{\pi+i\beta}{r}} & 0 \\[5mm] 
      0 & \displaystyle \frac{\sin\frac{\pi}{r}}
      {\sin\frac{\pi+i\beta}{r}} &  
      \displaystyle \frac{\sin\frac{i\beta}{r}}
      {\sin\frac{\pi+i\beta}{r}} & 0 \\[5mm] 
      0 & 0 & 0 & 1 \\
    \end{array} 
  \right) \;,
  \label{eq:dyra6}
\end{equation}
where the normalisation factor is
\begin{equation}
  \rho_r(\beta) = \frac{S_{2}^2(1+\frac{i\beta}{\pi} \;\vert\; r,2)}
  {S_{2}(\frac{i\beta}{\pi} \;\vert\; r,2) \,
    S_{2}(2+\frac{i\beta}{\pi} \;\vert\; r,2)} \;.
  \label{eq:normdyrs}
\end{equation}
$S_{2}(x \vert \omega_{1},\omega_{2})$ is Barnes' double sine 
function of periods $\omega_{1}$ and $\omega_{2}$ defined by:
\begin{equation}
  S_{2}(x \vert \omega_{1},\omega_{2}) = 
  \frac{\Gamma_{2}(\omega_{1}+\omega_{2}-x \;\vert\; \omega_{1},\omega_{2})} 
  {\Gamma_{2}(x \;\vert\; \omega_{1},\omega_{2})} \;,
\end{equation}
where $\Gamma_r$ is the multiple Gamma function of order $r$
given by
\begin{equation}
  \Gamma_{r}(x \vert \omega_{1},\dots,\omega_{r}) = \exp
  \left(
    \frac{\partial}{\partial s} \sum_{n_{1},\dots,n_{r} \ge 0}
    (x+n_{1}\omega_{1}+\dots+n_{r}\omega_{r})^{-s}\Bigg\vert_{s=0}
  \right) \;.
\end{equation}
The algebra $\dy{r}{V6}$ is defined by
\begin{equation}
  R_{12}(\beta_{1}-\beta_{2}) \, L_1(\beta_{1}) \, L_2(\beta_{2}) = 
  L_2(\beta_{2}) \, L_1(\beta_{1}) \, 
  R_{12}^{*}(\beta_{1}-\beta_{2}) \,,
  \label{eq232}
\end{equation}
where $R^{*}_{12}(\beta,r) \equiv R_{12}(\beta,r-c)$.  \\
The Lax matrix $L$ must now be expanded on both positive \emph{and}
negative powers as
\begin{equation}
L(\beta) = \sum_{k \in \ZZ} L_k \beta^{-k} \;.
\end{equation}
A presentation similar to the double Yangian case is achieved by
introducing the following two matrices: 
\begin{eqnarray}
  \label{eq235}
  && L^+(\beta) \equiv L(\beta - \,i\pi c) \,, \\
  && L^-(\beta) \equiv \sigma_3 \, L(\beta-i\pi r) \, \sigma_3
  \,.
\end{eqnarray}
They obey coupled exchange relations following from (\ref{eq232}) and
periodicity/unitarity properties of the matrices
$R_{12}$ and $R^{*}_{12}$:
\begin{eqnarray}
  && R_{12}(\beta_{1}-\beta_{2}) \, L^\pm_1(\beta_{1}) \,
  L^\pm_2(\beta_{2}) = L^\pm_2(\beta_{2}) \, L^\pm_1(\beta_{1}) \,
  R^{*}_{12}(\beta_{1}-\beta_{2}) \,, \\
  && R_{12}(\beta_{1}-\beta_{2}-\,i\pi c) \,
  L^+_1(\beta_{1}) \, L^-_2(\beta_{2}) = L^-_2(\beta_{2}) \, L^+_1(\beta_{1})
  \, R^{*}_{12}(\beta_{1}-\beta_{2}) \,.
  \label{eq236}
\end{eqnarray}
Contrary to the case of the double Yangian, the matrices $L^+$ and
$L^-$ are \emph{not} independent. Note also that, due to conflicting 
conventions, the $r\rightarrow\infty$ limit of
$L^\pm$ in $\dy{r}{V6}$ corresponds to $L^\mp$ in $\dy{}{}$.

\subsection{Deformed double Yangian $\dy{r}{V8}$}
The $R$-matrix of the deformed double Yangian $\dy{r}{V8}$ is obtained
as the scaling limit of the $R$-matrix of the elliptic algebra $\elpa$
\cite{Konno,JKM}. It reads
\begin{equation}
  R_{V8}(\beta,r) = \rho_r(\beta) \left( \begin{array}{cccc}
  \displaystyle  
  \frac{\cos\frac{i\beta}{2r} \; \cos\frac{\pi}{2r}}
  {\cos\frac{\pi+i\beta}{2r}}  
  & 0 & 0 & \displaystyle -\frac{\sin\frac{i\beta}{2r}\;\sin\frac{\pi}{2r}} 
  {\cos\frac{\pi+i\beta}{2r}} \\
  0 & \displaystyle \frac{\sin\frac{i\beta}{2r}\;\cos\frac{\pi}{2r}} 
  {\sin\frac{\pi+i\beta}{2r}} & \displaystyle \frac{\cos\frac{i\beta}{2r} \; 
    \sin\frac{\pi}{2r}} {\sin\frac{\pi+i\beta}{2r}} & 0 \\
  0 & \displaystyle \frac{\cos\frac{i\beta}{2r}\;\sin\frac{\pi}{2r}} 
  {\sin\frac{\pi+i\beta}{2r}} & \displaystyle \frac{\sin\frac{i\beta}{2r} \; 
    \cos\frac{\pi}{2r}} {\sin\frac{\pi+i\beta}{2r}} & 0 \\
  \displaystyle -\frac{\sin\frac{i\beta}{2r} \; \sin\frac{\pi}{2r}} 
  {\cos\frac{\pi+i\beta}{2r}} & 0 & 0 & \displaystyle 
  \frac{\cos\frac{i\beta}{2r}\;\cos\frac{\pi}{2r}}
  {\cos\frac{\pi+i\beta}{2r}} \\  
\end{array} \right) \;.
\label{eq:dyra8}
\end{equation}
It is also obtained from the $R$-matrix of $\dy{r}{V6}$ by a gauge
transformation \cite{Konno}. The algebra $\dy{r}{V8}$ is defined by
the same relation as $\dy{r}{V6}$, albeit with the matrix $R_{V8}$,
and the same type of Lax matrix with positive and negative modes. 

\subsection{Deformed double Yangian $\dy{r}{F}$}
The $R$-matrix of $\dy{r}{F}$ is given by
\begin{equation}
  R(\beta;r) = \rho_r(\beta) 
  \left( 
    \begin{array}{cccc}  
      1 & 0 & 0 & 0 \\
      0 &
      \displaystyle\frac{\sin\frac{i\beta}{r}}{\sin\frac{\pi+i\beta}{r}} &  
      e^{\beta/r}
      \displaystyle\frac{\sin\frac{\pi}{r}}{\sin\frac{\pi+i\beta}{r}}
      &  
      0 \\
      0 & \quad e^{-\beta/r} 
      \displaystyle\frac{\sin\frac{\pi}{r}}{\sin\frac{\pi+i\beta}{r}}
      \quad & \quad  
      \displaystyle\frac{\sin\frac{i\beta}{r}}{\sin\frac{\pi+i\beta}{r}} 
      \quad & 0 
      \\
      0 & 0 & 0 & 1 \\
    \end{array} 
  \right) \;.
  \label{eq:dyrF}
\end{equation}
The normalisation factor is the same as for $\dy{r}{V6}$. The
definition of the algebra and the Lax operator are again formally identical.

\section{Twist from $\dy{}{}$ to $\dy{r}{}$:
  representation formula}
\setcounter{equation}{0}

\subsection{A notation for $P_{12}$ invariant matrices}
Let us denote by $M(b^+,b^-)$ the $4\!\times\!4$ matrix given by
\begin{equation}
  M(b^+,b^-) \equiv
  \left(
    \begin{array}{cccc}
      1&0&0&0\\
      0& \frac12(b^++b^-) & \frac12(b^+-b^-) &0\\[2mm]
      0& \frac12(b^+-b^-) & \frac12(b^++b^-) &0\\
      0&0&0&1\\
    \end{array}
  \right) \;.
  \label{eq:Mbpm}
\end{equation}
With this definition, we have $M(a,b) M(a',b')= M(aa',bb')$ and
$M(a,b)^{-1} = M(a^{-1},b^{-1})$.

\medskip
\noindent
Now,
\begin{equation}
  R[\dy{}{}](\beta) = \rho(\beta) M
  \left(1,\frac{i \beta -\pi}{i \beta +\pi}
\right) \;.
\label{eq:3.2}
\end{equation}

\medskip
\noindent
We have
$R[\dy{r}{V6}](\beta) = \rho_r(\beta) M(b_r^+,b_r^-)$, with
\begin{eqnarray}
  b_r^+ &=& \frac{\cos\frac{i \beta -\pi}{2r}}
  {\cos\frac{i \beta +\pi}{2r}}
  = \frac{\Gamma_1(r + \frac{i \beta }{\pi}+1 | 2r)
    \Gamma_1(r - \frac{i \beta }{\pi}- 1 | 2r)}
  {\Gamma_1(r + \frac{i \beta }{\pi} -1 | 2r)
    \Gamma_1(r - \frac{i \beta }{\pi}+ 1 | 2r)} \;,
  \\[3mm]
  b_r^- &=& \frac{\sin\frac{i \beta -\pi}{2r}}
  {\sin\frac{i \beta +\pi}{2r}}
  = \frac{\Gamma_1(\frac{i \beta }{\pi}+1 | 2r)
    \Gamma_1(2r - \frac{i \beta }{\pi}- 1 | 2r)}
  {\Gamma_1(\frac{i \beta }{\pi} -1 | 2r)
    \Gamma_1(2r - \frac{i \beta }{\pi}+ 1 | 2r)}
  \\
  &=&
  \frac{\Gamma_1(2r+ \frac{i \beta }{\pi}+1 | 2r)
    \Gamma_1(2r - \frac{i \beta }{\pi}- 1 | 2r)}
  {\Gamma_1(2r + \frac{i \beta }{\pi} -1 | 2r)
    \Gamma_1(2r - \frac{i \beta }{\pi}+ 1 | 2r)}
  \cdot \frac{i \beta -\pi}{i \beta +\pi} \ .
\end{eqnarray}

\subsection{The linear equation in representation}

We remark that the normalisation factor of $\dy{r}{V6}$ can be
rewritten as:
\begin{equation}
  \rho_r(\beta) 
  = \rho_F(-\beta;r) \rho(\beta) \rho_F(\beta;r)^{-1}
\label{eq:3.6}
\end{equation}
with
\begin{equation}
  \rho_F(\beta) 
  = \frac{\Gamma_{2}(\frac{i\beta}{\pi} + 1 + r \;\vert\; 2,r)^2}
  {\Gamma_{2}(\frac{i\beta}{\pi} + r \;\vert\; 2,r)
    \Gamma_{2}(\frac{i\beta}{\pi} + 2 + r \;\vert\; 2,r)}
\;.
\end{equation}

\noindent
Equations (\ref{eq:3.2}-\ref{eq:3.6}) allow us to write:
\begin{equation}
  R[\dy{r}{V6}] = F_{21}(-\beta) R[\dy{}{}]
  F_{12}(\beta)^{-1} \; .
\end{equation}
Using the notation (\ref{eq:Mbpm}), $F_{12}(\beta)$ is given by
\begin{equation}
  F_{12}(\beta) = \rho_F(\beta) \cdot M
  \left(
    \frac{\Gamma_1(\frac{i \beta }{\pi} + r - 1 | 2r)}
    {\Gamma_1(\frac{i \beta }{\pi} + r + 1 | 2r)},
    \frac{\Gamma_1(\frac{i \beta }{\pi} + 2r - 1 | 2r)}
    {\Gamma_1(\frac{i \beta }{\pi} + 2r + 1 | 2r)}
  \right)\;.
\end{equation}

\noindent
This twist-like matrix reads 
\begin{eqnarray}
  F_{12}(\beta) &=& \rho_F(\beta) \prod_{n=1}^\infty
  M
  \left(
    1,\frac{i \beta +\pi+2n\pi r}{i \beta -\pi+2n\pi r}
  \right)
  M
  \left(
    \frac{i \beta +\pi+(2n-1)\pi r}{i \beta -\pi+(2n-1)\pi r},1
  \right)
  \\
  &=&
  \prod_{n=1}^\infty
  R(\beta-i(2n)\pi r)^{-1} \ \tau(  R(\beta-i(2n-1)\pi r)^{-1} )
\end{eqnarray}
with
\begin{equation}
  \tau(M(a,b))=M(b,a)
\end{equation}
and where, unless differently specified, $R$ is the
$R$-matrix of $\dy{}{}$.
One uses here the representation of $\rho_F(\beta)$ as an infinite product
\begin{equation}
  \rho_F(\beta) = \prod_{n=1}^\infty  \rho(\beta - i n\pi r)^{-1}  \;.
\end{equation}
\noindent
The automorphism $\tau$ may be interpreted as the adjoint action of
$(-1)^{\frac12 h_0^{(1)}}$, so that
\begin{eqnarray}
  F_{12}(\beta) &=& \prod_{n=1}^\infty
  R(\beta-i(2n)\pi r)^{-1}
  \Ad\left((-1)^{\frac12 h_0^{(1)}}\right) R(\beta-i(2n-1)\pi r)^{-1}
  \\
  &=&
  \prod_{n=1}^\infty 
  \Ad\left((-1)^{\frac{n}2 h_0^{(1)}}\right) R(\beta-in\pi r)^{-1} \,
  .
  \label{eq:pi}
\end{eqnarray}
Hence $F$ is solution of the difference equation
\begin{equation}
  F(\beta-i\pi r) = (-1)^{-\frac12 h_0^{(1)}} F(\beta) (-1)^{\frac12 h_0^{(1)}}
  \cdot R(\beta-i\pi r) \;.
  \label{eq:label_rouge}
\end{equation}

\noindent
It would be tempting to relate the automorphism $\tau$ to the one used in
\cite{JKOS}, although the naive scaling of the latter does not give
back the former. For instance our $\tau$ is inner not outer. 
\\
All the infinite products are logarithmically divergent. They are
consistently regularised by the $\Gamma_1$ and $\Gamma_2$
functions. In particular, $\lim\limits_{r\rightarrow \infty}F = M(1,1)
= \un_4$.

\section{The universal form of  $\cF[\dy{}{} ; \dy{r}{V6}]$}
\setcounter{equation}{0}

We construct a universal twist $\cF$ from $\dy{}{}$ to
$\dy{r}{V6}$,  such that 
\begin{equation}
  F(\beta_1-\beta_2) = \pi_{\beta_1}\otimes\pi_{\beta_2} (\cF) \;.
\end{equation}
%% Define the spectral parameter dependent twist and
%% $R$-matrix  as 
%% \begin{eqnarray}
%%   \cF(\beta) &=& \Ad(e^{\frac{i\beta}{\pi} d} \otimes \un) (\cF)\\
%%   \cR(\beta) &=& \Ad(e^{\frac{i\beta}{\pi} d} \otimes \un) (\cR)
%% \end{eqnarray}
The form of the difference equation (\ref{eq:label_rouge}) 
obeyed by the conjectural
representation of the twist, together with the
known generic structures of linear equations obeyed by universal
twists \cite{BR,ABRR,JKOS} lead us to postulate the following linear
equation for $\cF$:
\begin{equation}
  \label{eq:linequni}
  \cF \equiv \cF(r) = \Ad(\phi^{-1}\otimes \un) (\cF) \cdot \cC
\end{equation}
with 
\begin{eqnarray}
  \phi &=& (-1)^{\frac12 h_0} e^{(r+c)d} \;, \\
  \cC &\equiv& e^{-\alpha c\otimes d - \gamma d\otimes c} \cR \;.
\end{eqnarray}
%% implying for $\cF(\beta)$
%% \begin{equation}
%%   \cF(\beta) = (-1)^{-\frac12 h_0^{(1)}} \cF(\beta+i\pi (r+c_1)) 
%%   (-1)^{\frac12 h_0^{(1)}}
%%   \cdot e^{-\alpha c\otimes d - \gamma d\otimes c}  \cR(\beta) \;.
%% \end{equation}
We now prove the consistency of these postulates. 
We will use the following preliminary properties:
\begin{itemize}
\item The operator $d$ in the double Yangian
  $\dy{}{}$ is defined by  $[d,e(u)]=\frac{d}{du}e(u)$ (see
  \cite{Kh}). The evaluation representations are related through 
  $\pi_{\beta+\beta'} = \pi_\beta \circ
  \Ad(\exp(\frac{i\beta'}{\pi}d))$. 
\item The operator $d$ satisfies
  $\Delta(d)=d\otimes 1 + 1\otimes d$.
%%\item The evaluation representation $\pi(\beta)$ with evaluation
%%  parameter $\beta$ is obtained as $\exp(-\beta d) \pi(0)$.
%%\item The shift of the spectral parameter $\beta$ by an amount
%%  $\beta'$ may hence be obtained as the action of 
%%  $\Ad(\exp(\frac{i\beta'}{\pi} d))$.
\item The generator $h_0$ of $\dy{}{}$ is such that
  \begin{equation}
    h_0 e(u) = e(u)(h_0+2)\;, \qquad  h_0 f(u) = f(u)(h_0-2) \;,
    \qquad [h_0, h(u)]=0 \;,
  \end{equation}
  and hence $\tau= Ad\left((-1)^{\frac12 h_0^{(1)}}\right)$ satisfies
  $\tau^2=1$.
\end{itemize}
The equation (\ref{eq:linequni}) can be solved by 
\begin{equation}\label{}
  \cF(r) = \prod_k^{\longleftarrow} \cF_k(r) \;,
  \qquad \cF_k(r) = \phi_1^k \cC_{12}^{-1}\phi_1^{-k} \;.
\end{equation}
It is easily seen that equation (\ref{eq:pi}) is the evaluation
representation of this universal formula.
\\
As in \cite{JKOS}, $\cF_k$ satisfy the following properties:
\begin{eqnarray}\label{}
  (\Delta\otimes\id)(\cF_k(r)) &=& \cF_k^{(23)}(r+c_1)
  \cF_k^{(13)}\left(r+c_2+\frac\alpha{k}c_2\right) \;,
  \\
  (\id\otimes\Delta)(\cF_k(r)) &=& \cF_k^{(12)}(r)
  \cF_k^{(13)}\left(r-\frac\gamma{k}c_2\right) \;,
\end{eqnarray}
and
\begin{equation}
  \cF_k^{(12)}(r)
  \cF_{k+l}^{(13)}\left(r+\frac{l-\gamma}{k+l}c_2\right)
  \cF_l^{(23)}(r+c_1)
  =
  \cF_l^{(23)}(r+c_1)
  \cF_{k+l}^{(13)}\left(r+\frac{l+\alpha}{k+l}c_2\right)
  \cF_k^{(12)}(r) \;.
\end{equation}

\noindent
It is then straightforward to follow \cite{JKOS} to prove the shifted cocycle
relation, provided that
$\alpha+\gamma = -1$.

\noindent
We then have 
\begin{equation}
  \cF^{(12)}(r) (\Delta\otimes\id)(\cF(r))
  = \cF^{(23)}\left(r+c^{(1)}\right) (\id\otimes\Delta) (\cF(r)) \;.
\end{equation}
It follows that ${\cR_{\dy{r}{V6}}}_{12} = \cF_{21} \cR_{12}
\cF_{12}^{-1}$  satisfies a shifted Yang--Baxter equation 
\begin{equation}
  \cR_{12}(r+c^{(3)}) \, 
  \cR_{13}(r) \,
  \cR_{23}(r+c^{(1)}) = 
  \cR_{23}(r) \,
  \cR_{13}(r+c^{(2)}) \, 
  \cR_{12}(r) \;,
  \label{eq:sYBE}
\end{equation}
and that $\dy{r}{V6}$ is a quasi-Hopf algebra with 
$\Delta^\cF(x) = \cF \Delta(x) \cF^{-1}$ and 
$\Phi_{123} = \cF_{23}(r) \cF_{23}(r+c^{(1)})^{-1}$.

\section{Twist to $\dy{r}{V8}$}
\setcounter{equation}{0}

\subsection{In representation}
The $R$-matrix of $\dy{r}{V6}$ and $\dy{r}{V8}$ are related by 
\begin{equation}
  R[\dy{r}{V8}] = K_{21} \; R[\dy{r}{V6}] \; K_{12}^{-1} \;,
  \label{le_bleue}
\end{equation}
where
\begin{equation}
  K = V \otimes V \qquad \mbox{with} \qquad V = \frac{1}{\sqrt{2}} 
  \left( 
    \begin{array}{rr} 1 & 1 \\ -1 & 1 \\ 
    \end{array} 
  \right) \;.
\end{equation}
This implies an isomorphism between $\dy{r}{V8}$ and $\dy{r}{V6}$ 
where the Lax operators are connected by $L_{V8} = V L_{V6} V^{-1}$.

\subsection{Universal form}

We identify $V$ with an evaluation representation of an element $g$ 
\begin{equation}
  V \equiv \pi_x(g)
\qquad
\mbox{with}
\qquad
  g = \exp\left(\frac{\pi}{2}(f_0-e_0)\right) \;.
  \label{eq:g}
\end{equation}
Since $e_0$ and $f_0$ lie in the undeformed Hopf subalgebra $sl(2)$ of
$\dy{}{}$ \cite{Kh}, the coproduct of $g$ reads 
\begin{equation}
  \Delta(g) = g\otimes g
  \label{eq:delta_g}
\end{equation}
so that
\begin{equation}
  g_1 g_2 \, \Delta^\cF(g^{-1}) \cF = g_1 g_2 \, \cF g_1^{-1} g_2^{-1}
  \;.
  \label{eq:gG}
\end{equation}
The two-cocycle 
$g_1 g_2 \Delta^\cF(g^{-1})$ is a coboundary (with respect to the
coproduct $\Delta^\cF$). 
In representation, (\ref{eq:gG}) is equal to the scaling limit of the
represented twist from $\uq{q}$ to $\elpa$ \cite{JKOS,Cla}.
\\
Note that this case is similar to the gauge transformation used in 
\cite{ESS} although $g$ is not purely Cartan. 

\noindent
It follows that
\begin{equation}
  \cR[\dy{r}{V8}] \equiv  g_1 g_2 \, \Delta_{21}^\cF(g^{-1}) \; 
  \cR[\dy{r}{V6}] \;
  \Delta_{12}^\cF(g) \, g_1^{-1} g_2^{-1}  
\end{equation}
satisfies the shifted Yang--Baxter equation (\ref{eq:sYBE}).
\\
To recover (\ref{le_bleue}), use (\ref{eq:gG}) and remark that 
$\pi_x\otimes\pi_x (g\otimes g)$ commutes with $R[\dy{}{}]$. 

\section{Twist to $\dy{r}{F}$}
\setcounter{equation}{0}

\subsection{Twist in representation}

The $R$-matrices of $\dy{r}{V6}$ and $\dy{r}{F}$ are related by:
\begin{equation}
  R[\dy{r}{F}](\beta_1-\beta_2) = 
  K^{(6)}_{21}(\beta_2,\beta_1) \; R[\dy{r}{V6}](\beta_1-\beta_2) \; 
  (K^{(6)}_{12})^{-1}(\beta_1,\beta_2) \;,
\end{equation}
where
\begin{equation}
  K^{(6)}(\beta_1,\beta_2) = V'(\beta_1) \otimes V'(\beta_2) 
  \qquad \mbox{with} 
  \qquad V'(\beta) = 
  \left( 
    \begin{array}{rr} e^{\frac{\beta}{2r} } & 0 \\ 0 &
      e^{-\frac{\beta}{2r}} \\  
    \end{array} 
  \right) \;.
  \label{eq:K6}
\end{equation}

\subsection{Universal twist}

Again, one identifies $V'(\beta)$ as the evaluation representation of an
algebra element
\begin{equation}
  V'(\beta) = \pi_{\beta}\left(g' \right) \;,
\end{equation}
where
\begin{equation}
  g' = \exp\left(\frac{h_1}{2r}\right) \;.
\end{equation}
One then defines the following shifted coboundary  
\begin{equation}
  \cK_{12}(r) = g'(r) \otimes g'(r+c^{(1)}) \; 
  \Delta^\cF({g'}^{-1}) \;.
  \label{eq:cK12}
\end{equation}
It obeys a shifted cocycle condition 
\begin{equation}
  \cK_{12}(r) \; (\Delta^\cF \otimes \id) \cK(r) = \cK_{23}(r+c^{(1)})
  \; (\id  \otimes \Delta^{\cF'}) \cK(r) \;,
  \label{eq:cock}
\end{equation}
with $\cF'_{23}(r) = \cF_{23}(r+c^{(1)})$, 
as a consequence of 
\begin{equation}
  (\Delta^\cF \otimes \id) \; \Delta^\cF({g'}^{-1}) = (\id \otimes
  \Delta^{\cF'})  \; \Delta^\cF({g'}^{-1}) \;,
  \label{eq:cock2}
\end{equation}
which is the coassociativity property for $\Delta^\cF$.

\noindent
Finally
\begin{equation}
  \cR[\dy{r}{F}] \equiv  \cK_{21}(r)
  \cR[\dy{r}{V6}] \;
  \cK_{12}^{-1}(r)
  \label{non}
\end{equation}
satisfies the shifted Yang--Baxter equation (\ref{eq:sYBE}).
Moreover, (\ref{non})
together with (\ref{eq:cK12}) show that $\dy{r}{F}$ and $\dy{r}{V6}$ are the
same quasi-Hopf algebra.

\bigskip
\bigskip
\bigskip

\noindent \textbf{Acknowledgements}
\\
This work was supported in part by CNRS and EC network contract number 
FMRX-CT96-0012.  
\\
M.R. was supported by an EPSRC research grant no. GR/K 79437
and CNR-NATO fellowship. 
\\
D.A., L.F. and E.R. are most grateful to RIMS for hospitality.
We thank warmfully 
M.~Jimbo, H.~Konno, T.~Miwa and J.~Shiraishi for fruitful and
stimulating discussions.
\\
We are also indebted to S.~Pakuliak for his enlightening comments. 
\\
J.A.  wishes to thank the LAPTH for its kind hospitality.

\end{document}